\documentclass[12pt]{amsart}
\usepackage[utf8]{inputenc}
\usepackage[T1]{fontenc}

\usepackage[]{amsfonts}
\usepackage{amssymb}
\usepackage{amsmath}
\def\couleur(#1 #2 #3)
	{
\special{ps::[end]}}

\def\bx#1{\setbox1=\hbox{\kern3pt{#1}\kern3pt}			
 \dimen1=\ht1 \advance\dimen1 by 3pt \dimen2=\dp1 \advance\dimen2 by 3pt
 \setbox1=\hbox{\vrule height\dimen1 depth\dimen2\box1\vrule}%
 \setbox1=\vbox{\hrule\box1\hrule}%
 \advance\dimen1 by .4pt \ht1=\dimen1
 \advance\dimen2 by .4pt \dp1=\dimen2 \box1\relax}

\def\wbb#1{\kern#1em}
\def\vci{\vrule  width.02em height1.47ex depth-.0ex}		
\def\11{{\rm\wbb{.2}\vci\wbb{-.37}1}}

\def\Supp{\mathop{\rm Supp}\nolimits}
\def\underset#1#2{\mathrel{\mathop{\kern0pt #2}\limits_{#1}}}

\def\overset#1#2{\mathrel{\mathop{\kern0pt #2}\limits^{#1}}}

\newtheorem{Thrm}{Theorem}[section]
\newtheorem{Lmm}[Thrm]{Lemma}
\newtheorem{Dfnt}[Thrm]{Definition}
\newtheorem{Prps}[Thrm]{Proposition}
\newtheorem{Crll}[Thrm]{Corollary}
\newtheorem{Rmrq}[Thrm]{Remark}

\begin{document}

\title{On estimates for the $\displaystyle \bar \partial $ equation in Stein manifolds.}

\author{Eric Amar}

\address{UFR Math. Info. Universit\'e de Bordeaux. 351, Cours de la Lib\'eration, 33405, Talence France}

\date{}

\maketitle
 \ \par 
\ \par 
\ \par 

\tableofcontents
\ \par 
\ \par 
\renewcommand{\abstractname}{Abstract}

\begin{abstract}
We generalize to intersection of strictly $c$-convex domains
 in Stein manifold, $\displaystyle L^{r}-L^{s}$ and Lipschitz
 estimates for the solutions of the $\displaystyle \bar \partial
 $ equation obtained by Ma and Vassiliadou for domains in $\displaystyle
 {\mathbb{C}}^{n}.$ For this we use a Docquier-Grauert holomorphic
 retraction plus the raising steps method I introduce earlier.
 This gives results in the case of intersection of domains with
 low regularity, $\displaystyle {\mathcal{C}}^{3},$ for their boundary.\ \par 
\end{abstract}

\section{Introduction.}
\quad The solutions with   $\displaystyle L^{r}$  and Lipschitz estimates
 of the  equation $\displaystyle \bar \partial u=\omega ,\ \bar
 \partial \omega =0$   are known to be very important in complex
 analysis and geometry.\ \par 
\quad The first results of this kind were obtained by the use of solving
 kernels: Grauert-Lieb~\cite{GrauLieb70}, Henkin~\cite{Henkin70},
 Ovrelid~\cite{Ovrelid71}, Skoda~\cite{zeroSkoda}, Krantz~\cite{KrantzDbar76},
 in the case of  strictly pseudo-convex domains with $\displaystyle
 {\mathcal{C}}^{\infty }$ smooth boundary in $\displaystyle {\mathbb{C}}^{n},$
  with the exception of  Kerzman~\cite{Kerzman71} in the case
 of $\displaystyle (0,1)$  forms in  strictly pseudo-convex domains
 with $\displaystyle {\mathcal{C}}^{4}$ smooth boundary in Stein
 manifolds.\ \par 
\quad Here we shall be interested in strictly $c$-convex, s.c.c. for
 short, domain $D$ in a complex manifold. Such a domain is defined
 by a function  $\rho $ of class $\displaystyle {\mathcal{C}}^{3}$
 in a neighbourhood $U$ of $\bar D$  and such that  $\displaystyle
 i\partial \bar \partial \rho $  has at least $\displaystyle
 n-c+1$  strictly positive eigenvalues in $U.$\ \par 
\quad These domains in $\displaystyle {\mathbb{C}}^{n}$ have been studied
 in the case of smooth $\displaystyle {\mathcal{C}}^{\infty }$
 boundary by Fisher and Lieb~\cite{FischLieb74}.\ \par 
\quad Ma and Vassiliadou~\cite{MaVassiliadou00} obtained very nice
 estimates even in the case of intersections of s.c.c. domains
 with $\displaystyle {\mathcal{C}}^{3}$ boundary. I shall use
 their results here.\ \par 
\quad Quite recently C. Laurent-Thi\'ebaut~\cite{CLaurent13} got this
 kind of result for s.c.c. domains with smooth $\displaystyle
 {\mathcal{C}}^{\infty }$ boundary in complex manifold by use
 of the Grauert's method of "bumps".\ \par 
\quad Concerning the study of transverse intersection of domains, one
 can cite the works of Henkin and Leiterer~\cite{HenkLeit81},
 Menini~\cite{Menini97} for strictly pseudo convex domains and
 G. Schmalz~\cite{Schmalz89} and  Ma and Vassiliadou~\cite{MaVassiliadou00}
  for $c$-convex domains. C. Laurent-Thi\'ebaut and J. Leiterer~\cite{LeitLaur93}
 solved the $\displaystyle \bar \partial $ equation in a  case
 of intersection of s.c.c. domains more general than the one
 considered by Ma and Vassiliadou~\cite{MaVassiliadou00} but
 for bounded forms and they got solutions in Lipschitz spaces.
 It seems that the $\displaystyle L^{p}$ case is still open for
 their situation.\ \par 
\ \par 
\quad Let us state our first  result which is completely analogous
 to the one Ma and Vassiliadou~\cite{MaVassiliadou00} obtained
 for domains in $\displaystyle {\mathbb{C}}^{n}.$\ \par 

\begin{Thrm}
~\label{SE11}Let $\Omega $ be a Stein manifold of dimension $n.$
 Let $D$ be a strictly $c$-convex (s.c.c.) domain relatively
 compact with smooth $\displaystyle {\mathcal{C}}^{3}$ boundary
 in $\Omega .$ Let  $\omega $ be a $\displaystyle (p,q)$  form
 in $\displaystyle L^{r}_{p,q}(D),\ \bar \partial \omega =0,$
 with $\displaystyle 1<r<2n+2,\ c\leq q\leq n.$ Then there is
 a $\displaystyle (p,q-1)$  form $u$ in $\displaystyle L_{p,q-1}^{s}(D),$
 with $\displaystyle \ \frac{1}{s}=\frac{1}{r}-\frac{1}{2n+2},$
 such that $\displaystyle \bar \partial u=\omega .$\par 
\quad If  $\omega $ is in $\displaystyle L^{r}_{p,q}(D),\ \bar \partial
 \omega =0$  with  $\displaystyle r\geq 2n+2,\ c\leq q\leq n,$
 then there is a $\displaystyle (p,q-1)$  form $u$ in $\displaystyle
 \Lambda _{(p,q-1)}^{\epsilon }(\bar D)$ such that $\displaystyle
 \bar \partial u=\omega $ with $\displaystyle \epsilon =\frac{1}{2}-\frac{n+1}{r}.$
\end{Thrm}
The
 spaces $\displaystyle \Lambda _{(p,q-1)}^{\epsilon }(\bar D)$
 are the (isotropic) Lipschitz spaces of order $\epsilon $ and
 we set $\displaystyle \ \Lambda _{(p,q-1)}^{0}(\bar D):=L_{(p,q-1)}^{\infty
 }(D).$\ \par 
\quad It has to be noticed that the boundary regularity is just $\displaystyle
 {\mathcal{C}}^{3},$ so it seems that this is a new result in
 a Stein manifold for such a low regularity.\ \par 
\quad In the case of a $\displaystyle {\mathcal{C}}^{\infty }$ boundary
 regularity then this result is contained in C. Laurent-Thi\'ebaut~\cite{CLaurent13}
    corollary 2.11, but the proof here is completely different
 and, in some sense, "lighter" because it never uses Beals, Greiner
 and Stanton~\cite{BealsGreiStan87} heavy technology. We use
 for the analytic part kernels methods plus essentially geometric
 ones. Nevertheless we can recover the Sobolev estimates by a
 direct use of Beals, Greiner and Stanton~\cite{BealsGreiStan87}
 in the case of a $\displaystyle {\mathcal{C}}^{\infty }$ boundary
 regularity by theorem~\ref{SI23} here. This avoid the use of
 the "bumps method" but this is valid only in Stein manifolds
 although  C. Laurent-Thi\'ebaut~\cite{CLaurent13} results are
 valid in any complex manifold.\ \par 
\ \par 
\quad To state our next result, we need the definition of a $\displaystyle
 {\mathcal{C}}^{3}\ c$ convex intersection, still taken from~\cite{MaVassiliadou00}.\
 \par 

\begin{Dfnt}
~\label{SI14}A relatively compact domain $D$ in a Stein manifold
 $\displaystyle \Omega $ shall be called a $\displaystyle {\mathcal{C}}^{3}\
 c$-convex intersection if there exists a relatively compact
 neighbourhood $W$ in $\displaystyle \Omega $ of  $\bar D$ and
 a finite number of real $\displaystyle {\mathcal{C}}^{3}$ functions
 $\displaystyle \rho _{1},...,\ \rho _{N}$ where $\displaystyle
 n\geq N+3$ defined on $W$ such that $\displaystyle D=\lbrace
 z\in W::\rho _{1}(z)<0,...,\ \rho _{N}(z)<0\rbrace $ and the
 following are true:\par 
i) For $\displaystyle 1\leq i_{1}<\cdot \cdot \cdot <i_{l}\leq
 N$ the $1$-forms $\displaystyle d\rho _{i_{1}},...,\ d\rho _{i_{l}}$
 are $\displaystyle {\mathbb{R}}$-linearly independent on $\displaystyle
 \ \bigcap_{j=1}^{l}{\lbrace \rho _{i_{j}}\leq 0\rbrace .}$\par 
ii) For $\displaystyle 1\leq i_{1}<\cdot \cdot \cdot <i_{l}\leq
 N,$ for every $\displaystyle z\in \bigcap_{j=1}^{l}{\lbrace
 \rho _{i_{j}}\leq 0\rbrace },$ if we set $\displaystyle I:=(i_{1},...,\
 i_{l}),$ there exists a linear subspace $\displaystyle T_{z}^{I}$
 of $\displaystyle \Omega $ of complex dimension at least $\displaystyle
 n-c+1$ such that for $\displaystyle i\in I$ the Levi forms $\displaystyle
 L\rho _{i}$ restricted on $\displaystyle T_{z}^{I}$ are positive definite.
\end{Dfnt}
\quad We notice that, in $\displaystyle {\mathbb{C}}^{n},$ Ma and Vassiliadou
  need $\displaystyle N\leq n-2$ and here we need $\displaystyle
 N\leq n-3.$ Now we can state:\ \par 

\begin{Thrm}
~\label{SE12}Let $\Omega $ be a Stein manifold of dimension $n$
 and a  $\displaystyle {\mathcal{C}}^{3}\ c$-convex intersection
 $D$ such that $D$ is relatively compact in $\displaystyle \Omega
 .$ There exists a $\displaystyle \nu \in {\mathbb{N}}^{+}$ (which
 depends on the maximal number of non empty intersections of
 $\displaystyle \lbrace \rho _{j}=0\rbrace $) such that :\par 
if  $\omega $ is a $\displaystyle (p,q)$  form in $\displaystyle
 L^{r}_{p,q}(D),\ \bar \partial \omega =0$  with $\displaystyle
 q\geq c,\ 1<r<2n+2,$ then there is a $\displaystyle (p,q-1)$
  form $u$ in $\displaystyle L_{p,q-1}^{s}(D)$\!\!\!\! , such
 that $\displaystyle \bar \partial u=\omega $ with $\displaystyle
 \ \frac{1}{s}=\frac{1}{r}+\frac{1}{\lambda }-1,$ where $\displaystyle
 1\leq \lambda <\frac{2n+2\nu }{2n-1+2\nu }.$\par 
More precisely,\par 
i) For any $\displaystyle 1<r<2n+2\nu ,$ there exists $\displaystyle
 c_{r}(D),$ a positive constant such that\par 
\quad \quad \quad $\displaystyle \ {\left\Vert{u}\right\Vert}_{L^{s}_{(p,q-1)}(D)}\leq
 c_{r}(D){\left\Vert{\omega }\right\Vert}_{L^{r}_{(p,q)}(D)}$\par 
with $\displaystyle \ \frac{1}{s}=\frac{1}{r}-\frac{1}{2n+2\nu }.$\par 
ii) For $\displaystyle r\geq 2n+2\nu ,$ we have $\displaystyle
 \ {\left\Vert{u}\right\Vert}_{L^{\infty }_{(p,q-1)}(D)}\leq
 a_{r}(D){\left\Vert{\omega }\right\Vert}_{L^{r}_{(p,q)}(D)}$
 for some positive constant $\displaystyle a_{r}(D).$
\end{Thrm}
\quad This also seems to be new in case $\displaystyle {\mathbb{C}}^{n}$
 is replaced by a Stein manifold.\ \par 
\ \par 
\quad The results of Ma and Vassiliadou~\cite{MaVassiliadou00} give
  good estimates in case of domains in $\displaystyle {\mathbb{C}}^{n}.$
 The first point here was to pass from $\displaystyle {\mathbb{C}}^{n}$
 to a submanifold of $\displaystyle {\mathbb{C}}^{n}.$ To do
 this I was inspired by a nice paper of H. Rossi~\cite{Rossi76}
 on Docquier Grauert holomorphic retraction. The first result
 is based on it and is the following non optimal theorem.\ \par 

\begin{Thrm}
~\label{SE7}Let $M$ be a closed submanifold of dimension $d$
 of a Stein domain $\displaystyle U_{0}$ in $\displaystyle {\mathbb{C}}^{n}.$
  Let  $D$ be a s.c.c. domain relatively compact in $M\ (\bar
 D\subset M)$ with $\displaystyle {\mathcal{C}}^{3}$ boundary.
 Then, with $\displaystyle r\geq 2n+2,$ we can solve in $\displaystyle
 D$ the equation $\displaystyle \bar \partial u=\omega $ when
 $\displaystyle \bar \partial \omega =0$ and with $\displaystyle
 u\in \Lambda _{(p,q-1)}^{\epsilon }(\bar D)$ if $\displaystyle
 \omega \in L^{r}_{(p,q)}(D),\ c\leq q\leq n,$ with $\displaystyle
 \epsilon =\frac{1}{2}-\frac{n+1}{r}.$
\end{Thrm}
\quad Then we use the raising steps method~\cite{AmarSt13} (see also~\cite{HodgeCompact15}
 for more general operators than $\displaystyle \bar \partial
 $ and ~\cite{HodgeNonCompact15} in the non compact case). Let
 me recall it in this specific situation.\ \par 

\begin{Thrm}
~\label{SI21}Let $M$ be a closed complex manifold and $D$ a relatively
 compact domain in $\displaystyle M.$ Suppose there is $\displaystyle
 \delta >0$ and a finite covering $\displaystyle \lbrace U_{j}\rbrace
 _{j=1,...,N}$ of $\bar D$ such that :\par 
\quad (i) $\displaystyle \forall r>1,\ \forall \omega \in L^{r}_{(p,q)}(D),\
 \bar \partial \omega =0,\ \exists u_{j}\in L^{t}_{(p,q-1)}(D\cap
 U_{j})::\bar \partial u_{j}=\omega $ in $\displaystyle D\cap U_{j}$\par 
and $\displaystyle \ {\left\Vert{u_{j}}\right\Vert}_{L^{t}(D\cap
 U_{j})}\lesssim {\left\Vert{\omega }\right\Vert}_{L^{r}(D\cap
 U_{j})},$ with $\displaystyle \ \frac{1}{t}=\frac{1}{r}-\delta .$\par 
\quad (ii) $\displaystyle \exists s>1,\ \forall \omega \in L^{s}_{(p,q)}(\Omega
 ),\ \bar \partial \omega =0,\ \exists w\in L^{s}_{(p,q-1)}(\Omega
 )::\bar \partial w=\omega $ and $\displaystyle \ {\left\Vert{w}\right\Vert}_{L^{s}(\Omega
 )}\lesssim {\left\Vert{\omega }\right\Vert}_{L^{s}(\Omega )}.$\par 
Then there is a constant $\displaystyle c>0$ such that, for $\displaystyle
 r\leq s,$ if $\displaystyle \omega \in L^{r}_{(p,q)}(D),\ \bar
 \partial \omega =0,$ it exists $\displaystyle u\in L^{t}_{(p,q-1)}(D)$
 with $\displaystyle \lambda :=\min \ (\delta ,\frac{1}{r}-\frac{1}{s})$
 and $\displaystyle \ \frac{1}{t}=\frac{1}{r}-\lambda ,$ such
 that $\displaystyle \bar \partial u=\omega ,\ u\in L^{t}_{(p,q-1)}(D),\
 {\left\Vert{u}\right\Vert}_{L^{t}(D)}\leq c{\left\Vert{\omega
 }\right\Vert}_{L^{r}(D)}.$
\end{Thrm}
\quad The local estimates {\sl (i)} are given by "localizing s.c.c.
 domain", proposition~\ref{SE10}, plus the results of Ma and
 Vassiliadou, theorem~\ref{SE15} here. The global estimate {\sl
 (ii)}, the threshold, is given by the  $\displaystyle L^{r}-\Lambda
 ^{\epsilon }$ estimates done in theorem~\ref{SE7}. We get the
 same optimal results as for domains in $\displaystyle {\mathbb{C}}^{n}.$\ \par 

\begin{Thrm}
~\label{SE14}Let $M$ be a  complex submanifold of dimension $d$
 in $\displaystyle {\mathbb{C}}^{n}$ and a s.c.c. domain $D$
 such that $D$ is relatively compact with smooth boundary of
 class $\displaystyle {\mathcal{C}}^{3}$ in $M.$ Let $\omega
 $ be a $\displaystyle (p,q)$  form in $\displaystyle L^{r}_{p,q}(D),\
 \bar \partial \omega =0,\ c\leq q\leq n,$ with $\displaystyle
 1<r<2d+2.$ Then there is a $\displaystyle (p,q-1)$  form $u$
 in $\displaystyle L_{p,q-1}^{s}(D)$\!\!\!\! , with $\displaystyle
 \ \frac{1}{s}=\frac{1}{r}-\frac{1}{2d+2},$ such that $\displaystyle
 \bar \partial u=\omega .$\par 
\quad If $\displaystyle r\geq 2n+2$ then there is a $\displaystyle
 (p,q-1)$  form $u$ in $\displaystyle \Lambda _{(p,q-1)}^{\epsilon
 }(\bar D)$ such that $\displaystyle \bar \partial u=\omega $
 with $\displaystyle \epsilon =\frac{1}{2}-\frac{d+1}{r}.$
\end{Thrm}
\quad We follow exactly the same path to work with  $\displaystyle
 {\mathcal{C}}^{3}\ c$-convex intersection with again local estimates
 given by "localizing s.c.c. intersection", proposition~\ref{SI16}
 plus the results of Ma and Vassiliadou, theorem~\ref{SI17}.
 The global estimate (ii), the threshold, is given by the  $\displaystyle
 L^{r}-L^{\infty }$ estimates done in theorem~\ref{SI17} plus
 the generalization of a theorem of Rossi done in theorem~\ref{SI18}.\ \par 
\quad To pass to  Stein manifold, we use an embedding theorem  of Bishop
 and Narashiman (see theorem  5.3.9.  of H\"ormander~\cite{Hormander73})
 to see an abstract Stein manifold of dimension $d$ as a submanifold
 of $\displaystyle {\mathbb{C}}^{2d+1}.$ So we get our main results.\ \par 
\ \par 
This work will be presented in the following way.\ \par 
\quad $\displaystyle \bullet $ First we recall the estimates in the
 case of strictly $c$-convex domains in $\displaystyle {\mathbb{C}}^{n}$
 done by Ma and Vassiliadou~\cite{MaVassiliadou00}.\ \par 
\quad $\displaystyle \bullet $ We recall the Docquier Grauert holomorphic
 retraction on a complex submanifold $M$ of $\displaystyle {\mathbb{C}}^{n}.$\
 \par 
\quad $\displaystyle \bullet $ We extend a form $\omega $ from a domain
 $D$ s.c.c. in $M$ to a domain $E$ s.c.c. in $\displaystyle {\mathbb{C}}^{n}$
 by use of a generalization of a theorem of H. Rossi~\cite{Rossi76}.
 We then solve the form in $E$ by the known estimates in $\displaystyle
 {\mathbb{C}}^{n}.$ Then we show that the solution in $E$ can
 be restricted to $D$ to get a solution in $D$ with good enough
 estimates, for $\displaystyle r\geq 2n+2.$ This gives theorem~\ref{SE7}.\ \par 
\quad $\displaystyle \bullet $ We use the raising steps theorem with
 the threshold given by theorem~\ref{SE7}. So we have theorem~\ref{SE14}
 for the case of a submanifold of $\displaystyle {\mathbb{C}}^{n}.$\ \par 
\quad $\displaystyle \bullet $ Then by the same way, using ad-hoc 
 modifications of propositions done in the appendix, we get theorem~\ref{SE12}
 for s.c.c. intersections in the case of a submanifold of $\displaystyle
 {\mathbb{C}}^{n}.$\ \par 
\quad $\displaystyle \bullet $ By use of a theorem of Bishop and Narashiman,
 i.e. the proper embedding of a Stein manifold of dimension $d$
 in $\displaystyle {\mathbb{C}}^{2d+1},$ we get our main theorems~\ref{SE11},
 ~\ref{SE12} for any Stein manifold.\ \par 
\quad $\displaystyle \bullet $ Finally we prove technical results we
 need in the appendix.\ \par 
\ \par 
\quad I am indebted to C. Laurent-Thi\'ebaut who pointed to me the
 precise link between the work of Beals, Greiner and Stanton~\cite{BealsGreiStan87}
 and the existence of actual solutions for the $\displaystyle
 \bar \partial $ Neuman problem.\ \par 
\quad Moreover I thank the referee for his nice suggestions : in particular
 the proof of theorem~\ref{SE7} comes from a slight modification
 of one of them and it simplifies substantially my original proof.\ \par 

\section{Strictly $c$-convex domain in $\displaystyle {\mathbb{C}}^{n}.$}
\quad We shall use the nice estimates for a  smoothly $\displaystyle
 {\mathcal{C}}^{3}$ bounded  $c$ convex domains in $\displaystyle
 {\mathbb{C}}^{n}$  obtained by Ma and Vassiliadou~\cite{MaVassiliadou00}
 , lemma 5.3. in their paper.\ \par 

\begin{Thrm}
~\label{SE15}Let $D$ be a bounded s.c.c.  domain in $\displaystyle
 {\mathbb{C}}^{n}$ with a $\displaystyle {\mathcal{C}}^{3}$ defining
 function. Then\par 
\quad \quad \quad $\displaystyle \forall \omega \in L_{(p,q)}^{r},\ \bar \partial
 \omega =0,\ c\leq q\leq n,$\par 
there exists $\displaystyle u\in L_{(p,q-1)}^{s}(D),		\ \frac{1}{s}=\frac{1}{r}-\frac{1}{2n+2},$
  with the following properties:\par 
\quad i)  if $\displaystyle 1\leq r<2n+2,\ \bar \partial u=\omega $
  in the sense of currents in $D.$\par 
\quad ii)  if $\displaystyle r=1,\ u\in L_{(p,q-1)}^{\frac{2n+2}{2n+1}-\eta
 }$  for any $\displaystyle \eta >0.$\par 
\quad iii)  if $\displaystyle 2n+2\leq r\leq \infty ,\ u\in \Lambda
 _{(p,q-1)}^{\epsilon }(\bar D)$  with $\displaystyle \epsilon
 =\frac{1}{2}-\frac{n+1}{r}.$
\end{Thrm}
They also prove results in the case of intersections.\ \par 

\begin{Thrm}
~\label{SI17}Let a  $\displaystyle {\mathcal{C}}^{3}\ c$-convex
 intersection domain $D$ be such that $D$ is relatively compact
 in $\displaystyle {\mathbb{C}}^{n}.$ Then there exists a $\displaystyle
 \nu \in {\mathbb{N}}^{+}$ (which depends on the maximal number
 of non empty intersections of $\displaystyle \lbrace \rho _{j}=0\rbrace
 $) such that :\par 
if  $\omega $ a $\displaystyle (p,q)$  form in $\displaystyle
 L^{r}_{p,q}(D),\ \bar \partial \omega =0$  with $\displaystyle
 q\geq c,\ 1<r<2n+2\nu ,$  there is a $\displaystyle (p,q-1)$
  form $u$ in $\displaystyle L^{s}(D),$ such that $\displaystyle
 \bar \partial u=\omega $ with $\displaystyle \ \frac{1}{s}=\frac{1}{r}+\frac{1}{\lambda
 }-1,$ where $\displaystyle 1\leq \lambda <\frac{2n+2\nu }{2n-1+2\nu }.$\par 
More precisely\par 
i) For any $\displaystyle 1<r<2n+2\nu ,$ there exists $\displaystyle
 c_{r}(D)$ positive constant such that\par 
\quad \quad \quad $\displaystyle \ {\left\Vert{u}\right\Vert}_{L^{s}_{(p,q-1)}(D)}\leq
 c_{r}(D){\left\Vert{\omega }\right\Vert}_{L^{r}_{(p,q)}(D)}$\par 
with $\displaystyle \ \frac{1}{s}=\frac{1}{r}-\frac{1}{2n+2\nu }.$\par 
ii) For $\displaystyle r\geq 2n+2\nu ,$ we have $\displaystyle
 \ {\left\Vert{u}\right\Vert}_{L^{\infty }_{(p,q-1)}(D)}\leq
 a_{r}(D){\left\Vert{\omega }\right\Vert}_{L^{r}_{(p,q)}(D)}$
 for some positive constant $\displaystyle a_{r}(D).$
\end{Thrm}

\section{The Docquier - Grauert holomorphic retraction.}
\quad We have the Docquier-Grauert lemma~\cite{DocGrauert60} :\ \par 

\begin{Lmm}
Let $K$ be a compact subset of a closed complex submanifold $M$
 of dimension $d$ in $\displaystyle {\mathbb{C}}^{n}.$ There
 is a neighbourhood $U$ of $K$ and a holomorphic map $\displaystyle
 \pi \ :\ U\rightarrow U\cap M$ such that $\pi (\zeta )=\zeta
 $ for $\displaystyle \zeta \in U\cap M.$
\end{Lmm}
\quad In fact we have more (Rossi~\cite{Rossi76}, p 172) from the argument
 of Docquier-Grauert we have that the fibers $\pi ^{-1}\pi \zeta
 $ of $\pi $ intersect $M$ transversely at all points of $M$
 and are of dimension $\displaystyle n-d.$\ \par 
\ \par 
\quad Let $\displaystyle M$ be a complex submanifold of dimension $d$
 in $\displaystyle {\mathbb{C}}^{n}$ and $D$ a relatively compact
 domain strictly $c$-convex in $M.$ We have the following lemma.\ \par 

\begin{Lmm}
~\label{SI20}Let $\displaystyle \zeta \in \bar D,$ there is a
 neighborhood $U$ of $\zeta $ in $\displaystyle {\mathbb{C}}^{n}$
 and a bi-holomorphic mapping  $\displaystyle (U,\varphi ),\
 \varphi :U\rightarrow T:=\varphi (U),$ such that, with $\displaystyle
 z=(z_{1},...,z_{n})$ the coordinates in $T,$ we have: $\displaystyle
 \varphi (D)\bigcap T=\lbrace z_{d+1}=\cdot \cdot \cdot =z_{n}=0\rbrace
 $ and the retraction   $\displaystyle \tilde \pi :=\varphi \circ
 \pi \circ \varphi ^{-1}$ read in the mapping  $\varphi $ is
 given by $\displaystyle \tilde \pi (z)=(z_{1},...,\ z_{d},\
 0,...,0)$\!\!\!\! , i.e. this is the orthogonal projection onto
 the subspace of $\displaystyle z':=(z_{1},...,z_{d}).$  Moreover
 one can choose for $T$ a tube around $\displaystyle \varphi
 (M)$  of width  $\displaystyle \delta >0.$
\end{Lmm}
\quad Proof.\ \par 
The manifold $\displaystyle M$  is given, by use  of the retraction
  $\pi ,$ by the functions $\displaystyle f_{k}(\zeta ):=\zeta
 _{k}-\pi _{k}(\zeta ),\ k=1,...,n.$ We have if $\displaystyle
 \zeta \in M,\ \zeta -\pi (\zeta )=0$  ;  if $\displaystyle \zeta
 \notin M,\ \zeta -\pi (\zeta )\neq 0,$  because $\displaystyle
 \pi (\zeta )\in M.$  The transversality of the fibers with respect
 to $M$ at all points of $\bar D$ insures that the Jacobian of
 the map $\displaystyle f=(f_{1},...,f_{n})$  has rank $\displaystyle
 n-d,$  which is the complex co-dimension of $M.$ Take  a point
 $\displaystyle \zeta ^{0}\in \bar D,$ there are $\displaystyle
 n-d$  functions $\displaystyle f_{j}$  which are independent
 in a neighborhood  $U$ of $\displaystyle \zeta ^{0}.$  Re-numerating
 the functions $\displaystyle f_{j}$  and the variables $\displaystyle
 \zeta _{k},$  we may suppose that the  determinant $\displaystyle
 (\frac{\partial f_{j}}{\partial \zeta _{k}})_{j,k=d+1,...,n}$
  is different from zero.\ \par 
\quad Now we shall make the change of variables $\displaystyle z=\varphi
 (\zeta )$  with $\displaystyle z_{j}=\zeta _{j},\ j=1,\cdots
 ,d\ ;\ z_{j}=f_{j}(\zeta ),\ j=d+1,\cdots ,n.$  This is actually
 a change of variables because the Jacobian of $\varphi $ is
 different from zero in the open set $U.$ We have that  the mapping
  $\varphi $ is a bi-holomorphism from the open set $U$  onto
 the open set $\displaystyle T:=\varphi (U).$\ \par 
\quad Let $\displaystyle z'=(z_{1},\cdots ,z_{d})$  and $\displaystyle
 z"=(z_{d+1},\cdots ,z_{n})$ ; we have in $T$  that:\ \par 
\quad \quad \quad $\displaystyle N:=\varphi (M)=\lbrace z=(z',z")\in T::z"=0\rbrace .$\ \par 
Now take a tube around $N,\ T=\lbrace z=(z',z")::z"\in B((z',0),\
 \delta )\rbrace ,$ we call it again $T,$ and we still denote
 by $\displaystyle U$ the set $\displaystyle \varphi ^{-1}(T).$\ \par 
\quad We cover $\bar D$ by a finite number of these charts $\displaystyle
 (U_{j},\varphi _{j}).$ We note $\displaystyle N_{j}$ the manifold
 $\displaystyle N_{j}:=\varphi _{j}(M\cap U_{j})\subset T_{j}:=\varphi
 _{j}(U_{j})$ and, diminishing a little bit the $\displaystyle
 U_{j}$ if necessary, we can suppose that the width of the tubes
 $\displaystyle T_{j}$ around the $\displaystyle N_{j}$ is constant
 and equals $\displaystyle \delta >0.$ We know that there is
 a constant $\displaystyle \mu >0$ such that $\displaystyle \mu
 ^{-1}<J_{j}<\mu ,$ where $\displaystyle J_{j}$ is the Jacobian
 of $\displaystyle \varphi _{j},$ because there is a finite number
 of charts $\displaystyle (U_{j},\varphi _{j}).$ $\blacksquare $\ \par 
\quad In fact this lemma is a slightly more explicit local version
 of the proof of the Docquier-Grauert theorem in the book by
 Gunning and Rossi(~\cite{GunnRossi65}, Theorem 8, p.257) as
 was noticed by the referee.\ \par 
\quad We shall need this local version in order to get estimates in
 Lebesgue norms.\ \par 
\quad We denote $\displaystyle dV$ the Lebesgue measure on the manifold
 $\displaystyle M$ and by $\displaystyle dm$ the Lebesgue measure
 in $\displaystyle {\mathbb{C}}^{k}.$\ \par 
\quad We have the following basic lemma.\ \par 

\begin{Lmm}
~\label{SE1}Let  $f$ be a measurable function, positive on $M$\!\!\!\!
 , then\par 
\quad \quad \quad $\displaystyle \ \int_{U_{j}}{f\circ \pi (\zeta )dm(\zeta )}\leq
 \mu c(\delta )\int_{N_{j}}{f(z',0)dV(z')},$\par 
where $\displaystyle c(\delta ):=\left\vert{B(x,\delta )}\right\vert
 $ is the volume of the ball $\displaystyle B(x,\delta ).$
\end{Lmm}
\quad Proof.\ \par 
Let $f$ be a function in $\displaystyle L^{1}(U_{j})$  and $\tilde
 f$  this function read in the map $\displaystyle \varphi _{j},$
 i.e. $\displaystyle \tilde f:=f\circ \varphi _{j}^{-1},$ we get\ \par 
\quad \quad \quad \begin{equation}  \ \int_{U_{j}}{f(\zeta )dm(\zeta )}=\int_{N_{j}}{\lbrace
 \int_{B((z',0),\delta )}{\tilde f(z',z")J_{j}(z',z")dm(z")}\rbrace
 dV(z')}.\label{SI27}\end{equation}\ \par 
This is simply the change of variables formula because $\displaystyle
 \varphi _{j}(U_{j})=T_{j}=N_{j}{\times}B(\cdot ,\delta )$ and
 the Jacobian of $\displaystyle \varphi _{j}$ is $\displaystyle J_{j}.$ \ \par 
With the notation $\displaystyle z=(z',z"),\ z'$ the coordinates
 in $\displaystyle N_{j},\ z"$  the coordinates in the fibers,
 equation~(\ref{SI27}) gives:\ \par 
\quad \quad \quad $\displaystyle \ \int_{U_{j}}{f\circ \pi (\zeta )dm(\zeta )}=\int_{N_{j}}{\lbrace
 \int_{B(z',\delta )}{\tilde f(z',z")J_{j}(z',z")dm(z")}\rbrace dV(z')}.$\ \par 
Here we have $\displaystyle \tilde f(z',z")=f(z',0)$ because
 $\displaystyle \tilde \pi (z)=(z',0)$ hence the formula is now:\ \par 
\quad \quad \quad \quad \quad $\displaystyle \ \int_{U_{j}}{f\circ \pi (z)dm(z)}\leq \mu \int_{N_{j}}{f(z',0)\left\vert{B(z',\delta
 )}\right\vert dV(z')}=\mu c(\delta )\int_{N_{j}}{f(z',0)dV(z')}.$
 $\blacksquare $\ \par 
\quad We notice that the open set $\displaystyle U:=\bigcup_{j=1}^{N}{U_{j}}$
 contains $\bar D.$\ \par 

\section{Estimates in the case of a submanifold of $\displaystyle
 {\mathbb{C}}^{n}$ .}
\quad We shall show the following theorem:\ \par 

\begin{Thrm}
~\label{SI25}Let $M$ be a complex submanifold of dimension $d$
 in $\displaystyle {\mathbb{C}}^{n}$ and a s.c.c. domain $D$
 such that $D$  is relatively compact with smooth $\displaystyle
 {\mathcal{C}}^{3}$ boundary in $M.$ Let $\omega $ be a $\displaystyle
 (p,q)$  form in $\displaystyle L^{r}_{p,q}(D),\ \bar \partial
 \omega =0$  with  $\displaystyle r>2n+2,\ c\leq q\leq n.$ Then
 there is a $\displaystyle (p,q-1)$  form $u$ in $\displaystyle
 \Lambda _{(p,q-1)}^{\epsilon }(\bar D),\ \epsilon =\frac{1}{2}-\frac{n+1}{r},$
 such that $\displaystyle \bar \partial u=\omega .$
\end{Thrm}
\quad Proof.\ \par 
The idea is the following one:\ \par 
first we extend the form $\omega $ from $D$ to $E$ by $\displaystyle
 \tilde \omega =\pi ^{*}\omega .$ Then we solve the equation
 $\displaystyle \bar \partial \tilde u=\tilde \omega $ in $E$
 by theorem~\ref{SE15} with the estimates on $\displaystyle \tilde
 u.$ Then we restrict $\displaystyle \tilde u$ to $D$ to get
 the solution $\displaystyle u.$\ \par 
\quad Let us see that.\ \par 
The theorem~\ref{SE9} in the appendix, which generalizes to s.c.c.
 domains a theorem by Rossi~\cite{Rossi76} obtained for strictly
 pseudo convex domains, gives us the existence of a strictly
 $c$-convex domain $E\subset U$ in $\displaystyle {\mathbb{C}}^{n}$
 such that $\displaystyle \pi :\ \bar E\rightarrow \bar D.$ The
 open set $U$ here is $\displaystyle U:=\bigcup_{j=1}^{N}{U_{j}}$
 where $\displaystyle (U_{j},\varphi _{j})$ is the covering from
 lemma~\ref{SI20}. Now on we fix this s.c.c. domain $E.$\ \par 
\quad Let  $\omega $ be a $\displaystyle (p,q)$  form in $\displaystyle
 L_{p,q}^{r}(D),\ \bar \partial $  closed ; we extend it in $E$
 by use of the retraction $\pi $  in the following manner : $\displaystyle
 \tilde \omega :=\pi ^{*}\omega .$ With the notations of lemma~\ref{SI20}
 we start  by extending $\omega $ to $\displaystyle U\bigcap
 M$  by zero outside $\bar D$ ; the coefficients of $\tilde \omega
 $ can be written $\displaystyle f\circ \pi $  hence, applying
 lemma~\ref{SE1} to the functions  $\displaystyle 		\left\vert{f\circ
 \pi }\right\vert ^{r}$  we get $\displaystyle \ {\left\Vert{\tilde
 \omega }\right\Vert}_{L_{p,q}^{r}(U_{j})}\leq \mu c(\delta ){\left\Vert{\omega
 }\right\Vert}_{L_{p,q}^{r}(N_{j})}.$ We have only a finite number
 of open sets $\displaystyle U_{j}$ to cover $\displaystyle \bar
 D,$ so we get $\displaystyle \ {\left\Vert{\tilde \omega }\right\Vert}_{L_{p,q}^{r}(U)}\lesssim
 \mu c(\delta ){\left\Vert{\omega }\right\Vert}_{L_{p,q}^{r}(D)}.$
  Because $\displaystyle E\subset U$ we get $\displaystyle \
 {\left\Vert{\tilde \omega }\right\Vert}_{L_{p,q}^{r}(E)}\leq
 \mu c(\delta ){\left\Vert{\omega }\right\Vert}_{L_{p,q}^{r}(D)}.$\ \par 
Now with $\displaystyle r>2n+2,$ we can solve the equation $\displaystyle
 \bar \partial $ in the Lipschitz space by the theorem~\ref{SE15}, iii):\ \par 
\quad \quad \quad $\displaystyle \exists \tilde u\in \Lambda _{(p,q-1)}^{\epsilon
 }(\bar E)$  with $\displaystyle \epsilon =\frac{1}{2}-\frac{n+1}{r}$
 such that\ \par 
\quad \quad \quad \begin{equation}  \bar \partial \tilde u=\tilde \omega .\label{SI28}\end{equation}\
 \par 
Let $\displaystyle j\ :\ \bar D\rightarrow \bar E$ denote the
 inclusion map which is holomorphic. Notice that $\pi \circ j$
 is the identity map on $\displaystyle \bar D.$\ \par 
\quad Set $\displaystyle u:=j^{*}\tilde u.$ Then $\bar \partial u=j^{*}\bar
 \partial \tilde u=j^{*}\tilde \omega $ by equation~(\ref{SI28}). So\ \par 
\quad \quad \quad $\displaystyle \bar \partial u=j^{*}\tilde \omega =j^{*}\pi ^{*}\omega
 =\omega $\ \par 
because $\displaystyle j^{*}\pi ^{*}$ is the identity map on
 forms on $\displaystyle D.$\ \par 
Because $\displaystyle \tilde u$ has its coefficients in $\displaystyle
 \Lambda ^{\epsilon }(\bar E)$  then $\displaystyle u$ has its
 coefficients in $\displaystyle \Lambda _{(p,q-1)}^{\epsilon
 }(\bar D)$ because the restriction of a Lipschitz function is
 a Lipschitz function. $\blacksquare $\ \par 

\begin{Rmrq}
The theorem just proved with $\displaystyle r>2n+2,$ is enough
 to apply the raising steps method, and its proof is natural
 and simple. The case $\displaystyle \omega \in L_{p,q}^{2n+2}(D),\
 \bar \partial \omega =0,$ with a solution $u$ of the equation
 $\bar \partial u=\omega $ in $\displaystyle L_{p,q-1}^{\infty
 }(D)$ is not stated here, the problem being that the restriction
 to $D$ of a bounded function in $E$ is in general not even defined,
 opposite to the case of a Lipschitz function. Nevertheless this
 result is true and is a special case of theorem~\ref{SI34} done
 for s.c.c. intersection, with a more involved proof.
\end{Rmrq}
\quad Now we are in position to apply the raising steps theorem.\ \par 

\begin{Thrm}
~\label{SI24}Let $M$ be a complex submanifold of dimension $d$
 in $\displaystyle {\mathbb{C}}^{n}$and  $D$ be a s.c.c. domain
 which is relatively compact with smooth $\displaystyle {\mathcal{C}}^{3}$
 boundary in $M.$ Let $\omega $ be a $\displaystyle (p,q)$  form
  in $\displaystyle L^{r}_{p,q}(D),\ \bar \partial \omega =0$
 with $\displaystyle 1<r<2d+2,\ c\leq q\leq n.$ Then  there is
 a $\displaystyle (p,q-1)$  form $u$ in $\displaystyle L_{p,q-1}^{s}(D),$
 with $\displaystyle \ \frac{1}{s}=\frac{1}{r}-\frac{1}{2d+2},$
 such that $\displaystyle \bar \partial u=\omega .$
\end{Thrm}
\quad Proof.\ \par 
In order to have the local result for all points in $\bar D$
 we use the same method as in ~\cite{AmarSt13}, but with the
 proposition~\ref{SE10} and the results of Ma and Vassiliadou~\cite{MaVassiliadou00}.
 Let us see it.\ \par 
\quad Let $\zeta \in \partial D$ and $(V,\ \varphi )$ be a chart in
 a neighbourhood of $\zeta $ in $\displaystyle M$ and $\omega
 $ a $\displaystyle (p,q)$  form in $\displaystyle L^{r}_{p,q}(D),\
 \bar \partial \omega =0$ with $\displaystyle 1<r<2d+2,\ c\leq
 q\leq n.$ We read this situation in $\displaystyle {\mathbb{C}}^{d}$
 via the chart $\displaystyle (V,\ \varphi )$ so we have an open
 set $\displaystyle W:=\varphi (V)\subset {\mathbb{C}}^{d}$ and
 a piece of s.c.c. domain $\displaystyle \varphi (V\cap D)$ near
 the point $\displaystyle \eta :=\varphi (\zeta )\in {\mathbb{C}}^{d}$
 because the bi-holomorphic map $\varphi $ keeps the s.c.c. property.
 By use of the localizing proposition~\ref{SE10} there exist
 a s.c.c. domain $E\subset \varphi (V\cap D),$ with $\displaystyle
 {\mathcal{C}}^{3}$ boundary, which shares a part of its boundary
 near $\eta $ with the boundary of $\displaystyle \varphi (V\cap
 D).$ We read the form $\omega $ by $\varphi ,$ which gives us
 a $\displaystyle \tilde \omega :=\varphi ^{*}\omega $ still
 in $\displaystyle L^{r}_{p,q}(\varphi (V\cap D)),\ \bar \partial
 \tilde \omega =0$ hence $\displaystyle \tilde \omega \in L^{r}_{p,q}(E).$
 Now we apply the results of Ma and Vassiliadou~\cite{MaVassiliadou00},
 theorem~\ref{SE15} here, to get a  $\displaystyle (p,q-1)$ 
 form $\tilde u$ solution of the equation $\displaystyle \bar
 \partial \tilde u=\tilde \omega ,\ \tilde u\in L^{s}(E),\ $with
 $\displaystyle \ \frac{1}{s}=\frac{1}{r}-\frac{1}{2d+2}.$ Back
 to $D$ via $\displaystyle \varphi ^{-1}$ we have our local estimates
 : set $\displaystyle u:=(\varphi ^{-1})^{*}\tilde u,$ then $\displaystyle
 \bar \partial u=\omega ,\ u\in L_{p,q-1}^{s}(\varphi ^{-1}(E)).$\ \par 
\quad Hence we have the (i) of the raising steps theorem~\ref{SI21}.\ \par 
\ \par 
\quad We have the global result, i.e. the (ii) of the raising steps
 theorem~\ref{SI21} : set $\displaystyle t>2n+2\ ;$ if $\displaystyle
 \mu \in L^{t}_{p,q}(D),\ \bar \partial \mu =0$ then we have
 a solution $v$ in $\displaystyle \Lambda _{p,q}^{\epsilon }(\bar
 D)\subset L_{p,q}^{\infty }(D),$ with $\epsilon =t-(2n+2),$
 such that $\displaystyle \bar \partial v=\mu $  by use of theorem~\ref{SI25}.\
 \par 
\quad Now we take $\displaystyle \omega \in L^{r}_{p,q}(D),\ \bar \partial
 \omega =0,$ then we have that the optimal exponent for the solution
 $u$ of the equation $\bar \partial u=\omega $ is $\displaystyle
 s$ such that $\displaystyle \ \frac{1}{s}=\frac{1}{r}-\frac{1}{2d+2}\
 ;$  we choose any real $\displaystyle t$ such that $\displaystyle
 t>\max \ (2n+2,\ s)$ as a threshold  and, because  $\displaystyle
 L_{p,q-1}^{\infty }(D)\subset L^{t}_{(p,q-1)}(D),$  for $D$
 is a bounded domain, we  have a global solution to $\displaystyle
 \bar \partial v=\mu $ in $\displaystyle L^{t}_{(p,q-1)}(D)$
 if $\displaystyle \mu \in L^{t}_{(p,q)}(D).$ Now $\displaystyle
 s<t$ gives that $\displaystyle u\in L^{s}_{(p,q-1)}(D),$ by
 the raising steps theorem~\ref{SI21}, and this  ends the proof.
 $\displaystyle \blacksquare $\ \par 

\section{The case of  $\displaystyle {\mathcal{C}}^{3}\ c$-convex
 intersection.}
\quad We proceed exactly the same way than for just one s.c.c. domain.\ \par 
For the local estimates we use the localizing proposition~\ref{SI16}
 and we repeat the proof above. This is the point where we need
 to have at most $\displaystyle N=n-3$ domains in a Stein manifold
 of dimension $n,$ compare to $\displaystyle n-2$ domains in
 $\displaystyle {\mathbb{C}}^{n}.$ By use of Ma and Vassiliadou
 main theorem~\ref{SI17} and with $\displaystyle \nu \in {\mathbb{N}}^{+}$
 defined there, we get :\ \par 

\begin{Thrm}
~\label{SI33}Let $M$ be a complex submanifold of $\displaystyle
 {\mathbb{C}}^{n}$ of dimension $d$ and a  $\displaystyle {\mathcal{C}}^{3}\
 c$-convex intersection domain $D$ such that $D$ is relatively
 compact in $M.$ Let  $\omega $ be a $\displaystyle (p,q)$  form
 in $\displaystyle L^{r}_{p,q}(D),\ \bar \partial \omega =0$
  with $\displaystyle q\geq c,\ 1<r<2d+2\nu .$ Then there is
 finite covering $\displaystyle \lbrace U_{j}\rbrace _{j=1,...,N}$
 of $\bar D$ and $\displaystyle (p,q-1)$  form $u_{j}$ in $\displaystyle
 L^{s}(U_{j}\cap D)$\!\!\!\! , such that $\displaystyle \bar
 \partial u_{j}=\omega $ in $\displaystyle U_{j}\cap D$ with
 $\displaystyle \ \frac{1}{s}=\frac{1}{r}+\frac{1}{\lambda }-1$\!\!\!\!
 , where $\displaystyle 1\leq \lambda <\frac{2d+2\nu }{2d-1+2\nu }$ .\par 
More precisely,\par 
i) For any $\displaystyle 1<r<2d+2\nu ,$ there exists $\displaystyle
 c_{r}(D),$ a positive constant such that\par 
\quad \quad \quad $\displaystyle \ {\left\Vert{u_{j}}\right\Vert}_{L^{s}_{(p,q-1)}(U_{j}\cap
 D)}\leq c_{r}(D){\left\Vert{\omega }\right\Vert}_{L^{r}_{(p,q)}(U_{j}\cap
 D)}$\par 
with $\displaystyle \ \frac{1}{s}=\frac{1}{r}-\frac{1}{2d+2\nu }.$\par 
ii) For $\displaystyle r\geq 2d+2\nu ,$ we have $\displaystyle
 \ {\left\Vert{u_{j}}\right\Vert}_{L^{\infty }_{(p,q-1)}(U_{j}\cap
 D)}\leq a_{r}(D){\left\Vert{\omega }\right\Vert}_{L^{r}_{(p,q)}(U_{j}\cap
 D)}$ for some positive constant $\displaystyle a_{r}(D).$
\end{Thrm}
\ \par 
\quad Now for the global threshold, we copy the proof of theorem~\ref{SI25},
 replacing $\displaystyle r_{0}>\max \ (2n+2,\ s)$ by $\displaystyle
 r_{0}>\max \ (2n+2\nu ,\ s).$ We remark that we have not the
 Lipschitz estimates here but only $\displaystyle L^{\infty }$
 estimates and this is why we have to use the next lemma~\ref{c3}.\ \par 
\quad Let $E$ be  $\displaystyle {\mathcal{C}}^{3}\ c$-convex intersection
 domain in $\displaystyle {\mathbb{C}}^{n},\ E\subset U$ and,
 with $\displaystyle r\geq 2n+2,$ we can solve the $\displaystyle
 \bar \partial $  in the space $\displaystyle L^{\infty }(E)\
 :\ \bar \partial \tilde u=\tilde \omega ,\ \tilde u\in L_{p,q-1}^{\infty
 }(E)$ by the theorem~\ref{SI33}. Fix $\displaystyle \omega \in
 L^{r}_{(p,q)}(D),\ \bar \partial \omega =0,$ with $\displaystyle
 \ \tilde \omega $ as above, we have $\displaystyle \tilde u\in
 L_{p,q-1}^{\infty }(E)$ also fixed.\ \par 

\begin{Lmm}
~\label{c3} We have, with $\displaystyle j:\ D\rightarrow E$
 the canonical injection,\par 
\quad \quad \quad $\displaystyle \ {\left\Vert{j^{*}\tilde \omega -j^{*}\tilde
 \omega _{\epsilon }}\right\Vert}_{L^{r}(D_{\epsilon })}\underset{\epsilon
 \rightarrow 0}{\rightarrow }0.$
\end{Lmm}
\quad Proof.\ \par 
Recall that the coefficients $\displaystyle {\tilde a}_{I,J}$
 of $\displaystyle \tilde \omega $ verify $\displaystyle {\tilde
 a}_{I,J}(z)=a_{I,J}\circ \pi (z),$ if $\displaystyle a_{I,J}$
 is the corresponding coefficient of $\omega .$ We have, by definition
 of the convolution, noting $\displaystyle {\tilde a}_{I,J}^{\epsilon
 }$ a coefficient of $\displaystyle \tilde \omega _{\epsilon },$\ \par 
\quad \quad \quad $\displaystyle {\tilde a}_{I,J}^{\epsilon }(z):=\int_{{\mathbb{C}}^{n}}{a_{I,J}\circ
 \pi (z-\zeta )\chi _{\epsilon }(\zeta )dm(\zeta )},$\ \par 
so\ \par 
\quad \quad \quad $\displaystyle {\tilde a}_{I,J}^{\epsilon }(z)-{\tilde a}_{I,J}(z)=\int_{{\mathbb{C}}^{n}}{(a_{I,J}\circ
 \pi (z-\zeta )-a_{I,J}\circ \pi (z))\chi _{\epsilon }(\zeta
 )dm(\zeta )}.$\ \par 
Now $\displaystyle j^{*}$ is the operator of restriction to $D$
 so take $\displaystyle z\in D,$ then in a chart $\displaystyle
 (U_{j},\varphi _{j}),$ keeping the same notations for the functions
 read in this chart, with $\displaystyle z=(z',z''),\ \zeta =(\zeta
 ',\zeta ''),$\ \par 
\quad \quad \quad $\displaystyle {\tilde a}_{I,J}^{\epsilon }(z',0)-{\tilde a}_{I,J}(z',0)=\int_{{\mathbb{C}}^{n}}{(a_{I,J}(z'-\zeta
 ')-a_{I,J}(z'))\chi _{\epsilon }(\zeta )dm(\zeta )},$\ \par 
because here $\pi $ is the orthogonal projection on $\displaystyle
 (z',0),$ hence $\displaystyle {\tilde a}_{I,J}(z)=a_{I,J}(z').$\ \par 
So, decomposing the measure, we set\ \par 
\quad \quad \quad $\rho _{\epsilon }(\zeta '):=\int_{{\mathbb{C}}^{n-d}}{\chi _{\epsilon
 }(\zeta ',\zeta '')dm(\zeta '')},$\ \par 
and $\rho _{\epsilon }(\zeta ')$ is an approximate identity in
 $\displaystyle {\mathbb{C}}^{d}$ because\ \par 
\quad \quad \quad $\displaystyle \ {\left\Vert{\rho _{\epsilon }}\right\Vert}_{L^{1}({\mathbb{C}}^{d})}=\int_{{\mathbb{C}}^{n}}{\chi
 _{\epsilon }(\zeta ',\zeta '')dm(\zeta )}=1$\ \par 
and\ \par 
\quad \quad \quad $\displaystyle \forall \varphi \in {\mathcal{C}}({\mathbb{C}}^{d}),\
 \forall \zeta '\in {\mathbb{C}}^{d,\ }\varphi *\rho _{\epsilon
 }(\zeta ')\underset{\epsilon \rightarrow 0}{\rightarrow }\varphi
 (\zeta ').$\ \par 
So we have a convolution on $\displaystyle a_{I,J}(z')$ hence\ \par 
\quad \quad \quad $\displaystyle {\tilde a}_{I,J}^{\epsilon }(z',0)-{\tilde a}_{I,J}(z',0)=\int_{{\mathbb{C}}^{d}}{(a_{I,J}(z'-\zeta
 ')-a_{I,J}(z'))\rho _{\epsilon }(\zeta )dm(\zeta ')},$\ \par 
and, because the convolution is continuous on $\displaystyle
 L^{r},\ \forall r\geq 1,$\ \par 
\quad \quad \quad $\displaystyle \ {\left\Vert{{\tilde a}_{I,J}^{\epsilon }-{\tilde
 a}_{I,J}}\right\Vert}_{L^{r}(U_{j}\cap D)}\lesssim {\left\Vert{a_{I,J}^{\epsilon
 }-a_{I,J}}\right\Vert}_{L^{r}(U_{j}\cap D)}\underset{\epsilon
 \rightarrow 0}{\rightarrow }0.$\ \par 
Because we have only a finite number of charts $\displaystyle
 (U_{j},\varphi _{j})$ to cover $E,$ the proof is complete. $\blacksquare
 $\ \par 
\ \par 
\quad Now we get, still with $\displaystyle \nu \in {\mathbb{N}}^{+}$given
 by Ma and Vassiliadou main theorem~\ref{SI17}:\ \par 

\begin{Thrm}
~\label{SI34}Let $M$ be a complex submanifold of dimension $d$
 in $\displaystyle {\mathbb{C}}^{n}$ and  a  $\displaystyle {\mathcal{C}}^{3}\
 c$-convex intersection domain $D$ such that $D$ is relatively
 compact in $M.$ Let $\omega $ be a $\displaystyle (p,q)$  form
 in $\displaystyle L^{r}_{p,q}(D),\ \bar \partial \omega =0$
  with  $\displaystyle r\geq 2n+2\nu ,\ c\leq q\leq n.$ Then
 there is a $\displaystyle (p,q-1)$  form $v$ in $\displaystyle
 L^{\infty }_{(p,q-1)}(D)$ such that $\displaystyle \bar \partial v=\omega .$
\end{Thrm}
\quad Proof.\ \par 
The only difference with the proof of theorem~\ref{SI25} is that
 the restriction of a $\displaystyle L^{\infty }(E)$ function
 to $D$ is not even defined a priori.\ \par 
\quad So we regularise the solution $\displaystyle \tilde u\in L^{\infty
 }(E)$ given by Ma and Vassiliadou main theorem~\ref{SI17} by
 convolution with a smooth function $\chi $ such that $\displaystyle
 \chi (t)\in {\mathcal{C}}^{\infty }_{c}(\rbrack 0,1\lbrack ).$\ \par 
As usual we choose $\chi $ such that $\displaystyle \ \int_{{\mathbb{C}}^{n}}{\chi
 (\left\vert{z}\right\vert ^{2})dm(z)}=1$ and we set $\displaystyle
 \chi _{\epsilon }(z):=\frac{1}{\epsilon ^{2n}}\chi (\frac{\left\vert{z}\right\vert
 ^{2}}{\epsilon ^{2}})$ and $\displaystyle \tilde u_{\epsilon
 }(z):=(\tilde u*\chi _{\epsilon })(z),$ which means that the
 convolution is done on the coefficients.\ \par 
Now we have that $\displaystyle \tilde u_{\epsilon }\in {\mathcal{C}}_{p,q-1}^{\infty
 }(E_{\epsilon })$ where $\displaystyle E_{\epsilon }:=\lbrace
 z\in E::d(z,E^{c})>\epsilon \rbrace $ and, because $\bar \partial
 \tilde u=\tilde \omega ,$\ \par 
\quad \quad \quad \begin{equation}  \bar \partial \tilde u_{\epsilon }=\tilde \omega
 _{\epsilon }\label{c0}\end{equation}\ \par 
with $\displaystyle \tilde \omega _{\epsilon }:=(\tilde \omega
 *\chi _{\epsilon })(z).$\ \par 
Moreover we have\ \par 
\quad \quad \quad $\displaystyle \forall \epsilon >0,\ {\left\Vert{\tilde u_{\epsilon
 }}\right\Vert}_{L_{p,q-1}^{\infty }(E_{\epsilon })}\leq {\left\Vert{\tilde
 u}\right\Vert}_{L_{p,q-1}^{\infty }(E)}$\ \par 
and\ \par 
\quad \quad \quad \begin{equation}  \ {\left\Vert{\tilde \omega -\tilde \omega
 _{\epsilon }}\right\Vert}_{L_{p,q}^{r}(E_{\epsilon })}\underset{\epsilon
 \rightarrow 0}{\rightarrow }0.\label{c1}\end{equation}\ \par 
Let $\displaystyle j\ :\ \bar D\rightarrow \bar E$ denote the
 inclusion map which is holomorphic. Notice that $\pi \circ j$
 is the identity map on $\displaystyle \bar D.$\ \par 
\quad Set $\displaystyle u_{\epsilon }:=j^{*}\tilde u_{\epsilon }.$
 Then $\bar \partial u_{\epsilon }=j^{*}\bar \partial \tilde
 u_{\epsilon }=j^{*}\tilde \omega _{\epsilon }$ by equation~(\ref{c0}).\ \par 
Now by~(\ref{c1}) we have, by use of lemma~\ref{c3}, which is
 necessary because in general the restriction is {\sl not} a
 continuous operator on $\displaystyle L^{r},$\ \par 
\quad \quad \quad $\displaystyle \ {\left\Vert{j^{*}\tilde \omega -j^{*}\tilde
 \omega _{\epsilon }}\right\Vert}_{L^{r}(D_{\epsilon })}\underset{\epsilon
 \rightarrow 0}{\rightarrow }0.$\ \par 
So\ \par 
\quad \quad \quad \begin{equation}  \bar \partial u_{\epsilon }=j^{*}\tilde \omega
 _{\epsilon }=j^{*}(\tilde \omega _{\epsilon }-\tilde \omega
 )+j^{*}\pi ^{*}\omega =j^{*}(\tilde \omega _{\epsilon }-\tilde
 \omega )+\omega ,\label{c2}\end{equation}\ \par 
because $\displaystyle j^{*}\pi ^{*}$ is the identity map on
 on forms on $\displaystyle D.$\ \par 
\quad Take the sequence $\displaystyle \lbrace u_{1/k}\rbrace _{k\in
 {\mathbb{N}}},$ then, because $\displaystyle L^{\infty }(D)$
 is the dual of $\displaystyle L^{1}(D),$ there is a sub-sequence
 $\displaystyle \lbrace v_{k}\rbrace _{k\in {\mathbb{N}}}$ of
 $\displaystyle \lbrace u_{1/k}\rbrace _{k\in {\mathbb{N}}}$
 *-weakly converging, i.e. converging against $\displaystyle
 L_{n-p,n-q+1}^{1}(D)$ forms to $\displaystyle v\in L_{p,q-1}^{\infty
 }(D).$\ \par 
Let $\varphi \in {\mathcal{C}}^{\infty }_{n-p,n-q+1}(D)\subset
 L_{n-p,n-q+1}^{1}(D)$ with compact support in $\displaystyle
 D,$ i.e. a test form. Then for $k\geq k_{0}$ big enough we have
 that $\Supp \varphi \subset E_{1/k_{0}}$ and by~(\ref{c2}) we have:\ \par 
\quad \quad \quad $\displaystyle \forall k\geq k_{0},\ (-1)^{p+q-1}{\left\langle{v_{k},\bar
 \partial \varphi }\right\rangle}={\left\langle{\bar \partial
 v_{k},\varphi }\right\rangle}={\left\langle{j^{*}(\tilde \omega
 _{1/k}-\tilde \omega ),\varphi }\right\rangle}+{\left\langle{\omega
 ,\varphi }\right\rangle},$\ \par 
and this is well defined because $\displaystyle \Supp \varphi
 \subset E_{1/k_{0}}.$ Letting $\displaystyle k\rightarrow \infty
 $ we get\ \par 
\quad \quad \quad $\displaystyle (-1)^{p+q-1}{\left\langle{v_{k},\bar \partial
 \varphi }\right\rangle}\rightarrow {\left\langle{\omega ,\varphi
 }\right\rangle}.$\ \par 
But $\displaystyle \ {\left\langle{v_{k},\bar \partial \varphi
 }\right\rangle}\underset{k\rightarrow \infty }{\rightarrow }{\left\langle{v,\bar
 \partial \varphi }\right\rangle},$ because $\bar \partial \varphi
 \in L_{n-p,n-q+1}^{1}(D),$ so we get  $\displaystyle \ (-1)^{p+q-1}{\left\langle{v,\bar
 \partial \varphi }\right\rangle}={\left\langle{\omega ,\varphi
 }\right\rangle}$ which means that $\bar \partial v=\omega $
 in the distributions sense. $\blacksquare $\ \par 

\begin{Rmrq}
We have no such estimates in the case $\displaystyle r<2n+2$
 because the limit of mean values in balls of a function in $\displaystyle
 L^{s}$ in $E,$ which is the case of $v,$ is no longer in $\displaystyle
 L^{s}(D)$ for $\displaystyle s<\infty ,$ in general, as can be easily seen.
\end{Rmrq}
\quad So we can apply the raising steps theorem to get the analogous
 results in the case of  a  $\displaystyle {\mathcal{C}}^{3}\
 c$-convex intersection domain $D$ such that $D$ is relatively
 compact in $M,$ a closed complex submanifold of $\displaystyle
 {\mathbb{C}}^{n}.$\ \par 

\section{Estimates in the case of a Stein manifold.}
\quad We can apply a theorem of Bishop and Narashiman (see theorem
  5.3.9.  of H\"ormander~\cite{Hormander73}) which tells us that,
 if $\Omega $  is a Stein manifold of dimension $d$\!\!\!\! ,
 there is an element $\displaystyle f\in {\mathcal{H}}(\Omega
 )^{2d+1}$ which defines a  regular injective and proper map
 from  $\Omega $ in $\displaystyle {\mathbb{C}}^{2d+1}$ . Denote
 $\displaystyle M:=f(\Omega )$ ; if $\displaystyle D'$  is the
 strictly $c$-convex domain in $\Omega $\!\!\!\! , relatively
 compact in $\Omega ,$ then  its image $\displaystyle D=f(D')$
  is a strictly $c$-convex domain in $M$\!\!\!\! . We can apply
  theorems~\ref{SI25} and~\ref{SI24}.\ \par 
\quad Of course the same is true for  $\displaystyle {\mathcal{C}}^{3}\
 c$-convex intersection in Stein manifold $M$ and we get our
 main theorems.\ \par 
\quad We get an easy corollary of our main theorems, (see Ma and Vassiliadou~\cite{MaVassiliadou00},
 corollary 1.). We have, because $D$ is relatively compact, the
 estimate $\displaystyle L^{2}-L^{2},$ and this gives :\ \par 

\begin{Crll}
Let $\Omega $ be a Stein manifold of dimension $n$ and a strictly
 $c$-convex domain $D,$ or a $\displaystyle {\mathcal{C}}^{3}\
 c$-convex intersection, such that $D$ is relatively compact
 with smooth $\displaystyle {\mathcal{C}}^{3}$ boundary in $\Omega
 .$ Then, for $\displaystyle q\geq c,$ the operator $\displaystyle
 \bar \partial \ :\ L^{2}_{(p,q-1)}(D)\rightarrow L^{2}_{(p,q)}(D)$
 has closed range.
\end{Crll}
\quad Proof.\ \par 
This is fairly well known: for instance theorem 1.1.1 in~\cite{Hormander65}.
 $\displaystyle \blacksquare $\ \par 
\quad Because the $\displaystyle L^{2}$ norm of the canonical solution
 of $\displaystyle \bar \partial u=\omega $ is smaller than the
 solution we obtain, this implies the $\displaystyle L^{2}$ existence
 of the $\displaystyle \bar \partial $-Neumann operator on strictly
 $c$-convex domains. We also have that the strictly $c$-convex
 condition implies the $\displaystyle Z(q)$ condition of Beals,
 Greiner and Stanton~\cite{BealsGreiStan87} for $\displaystyle
 c\leq q\leq n,$ hence we get an automatic improvement of regularity
 in the case of a $\displaystyle {\mathcal{C}}^{\infty }$ smoothly
 bounded s.c.c. domain, by theorems 2 and 4 in Beals, Greiner
 and Stanton~\cite{BealsGreiStan87}.\ \par 

\begin{Thrm}
~\label{SI23}Let $\Omega $ be a Stein manifold of dimension $n$
 and a strictly $c$-convex (s.c.c.) domain $D$ such that $D$
 is relatively compact with smooth $\displaystyle {\mathcal{C}}^{\infty
 }$ boundary in $\Omega .$ Let $\displaystyle k\in {\mathbb{N}}$
 and  $\omega $ a $\displaystyle (p,q)$  form in $\displaystyle
 W^{k,r}(D),\ \bar \partial \omega =0$  with $\displaystyle 1<r<2n+2,\
 c\leq q\leq n.$ Then there is a $\displaystyle (p,q-1)$  form
 $u$ in $\displaystyle W^{k+1/2,r}(D),$ such that $\displaystyle
 \bar \partial u=\omega .$\par 
\quad If $\epsilon >0$ and $\omega $ is in $\displaystyle \Lambda ^{\epsilon
 }_{p,q}(D),\ \bar \partial \omega =0$  with  $\displaystyle
 c\leq q\leq n,$ then there is a $\displaystyle (p,q-1)$  form
 $u$ in $\displaystyle \Lambda _{(p,q-1)}^{\epsilon +1/2}(\bar
 D)$ such that $\displaystyle \bar \partial u=\omega .$
\end{Thrm}
Here we use the notation $\displaystyle W^{k,r}(D)$ for the Sobolev
 space of functions whose derivatives of order less than $k$
 are in $\displaystyle L^{r}.$\ \par 
\ \par 
\quad We  notice that there is no hypothesis here in the case $\displaystyle
 r>2$ on the  compactness of the support of the form $\omega
 ,$  in contrast to the previous results we had in~\cite{AmarSt13}.\ \par 

\section{Appendix.}
\ \par 

\begin{Lmm}
~\label{SE6}Let $\displaystyle A,\ B$ be two self adjoint matrices
  such that $A$ has at least  $\displaystyle n-c+1$  strictly
 positive eigenvalues and $B$  is positive. Then $\displaystyle
 A+B$  has at least $\displaystyle n-c+1$  strictly positive eigenvalues.
\end{Lmm}
\quad Proof.\ \par 
Let $E$ be the space generated by the eigenvectors associated
 to the strictly positive eigenvalues of $A.$ Then $E$ has dimension
 at least $\displaystyle n-c+1.$ Let $\displaystyle S:=A+B,$
 because $B$ is positive, we get\ \par 
\quad \quad \quad $\displaystyle \forall x\in E,{\left\langle{Sx,x}\right\rangle}={\left\langle{Ax,x}\right\rangle}+{\left\langle{Bx,x}\right\rangle}>0.$\
 \par 
Now let $\displaystyle e_{1},...,e_{k}$ be the eigenvectors associated
 to the negative eigenvalues of $\displaystyle S.$ We set $\displaystyle
 F=\mathrm{s}\mathrm{p}\mathrm{a}\mathrm{n}\lbrace e_{1},...,e_{k}\rbrace
 ,$ we have that $F$ is invariant by $S$ and we have $\displaystyle
 \forall x\in F,\ {\left\langle{Sx,\ x}\right\rangle}\leq 0.$
 If the space $\displaystyle G:=E\cap F$ is of non zero dimension,
 we get $\displaystyle \forall x\in G,\ x\neq 0,\ {\left\langle{Sx,x}\right\rangle}>0$
 and $\displaystyle \ {\left\langle{Sx,x}\right\rangle}\leq 0$
 so a contradiction. Hence $\displaystyle \mathrm{d}\mathrm{i}\mathrm{m}G=0$
 and $\displaystyle \mathrm{d}\mathrm{i}\mathrm{m}F\leq \mathrm{c}\mathrm{o}\mathrm{d}\mathrm{i}\mathrm{m}E=c-1,$
 which means that $S$ has a least $\displaystyle n-c+1$ strictly
 positive eigenvalues. $\displaystyle \blacksquare $\ \par 
\quad The next proposition generalizes the one in~\cite{AmCoh84}, proposition
 1.1, done for the pseudo convex case.\ \par 

\begin{Prps}
~\label{SE10}(Localizing s.c.c. domain)Let $D$ be a strictly
 $c$-convex domain with $\displaystyle {\mathcal{C}}^{3}$boundary
 in $\displaystyle {\mathbb{C}}^{n}.$ Let $\displaystyle \zeta
 \in \partial D,\ U$ a neighbourhood of $\zeta $ in $\displaystyle
 {\mathbb{C}}^{n}$ and $\displaystyle B(\zeta ,r)$ a ball centered
 at $\zeta $ and of radius $r$ such that $\displaystyle B(\zeta
 ,3r)\subset U$ ; then there is a domain $\displaystyle \ \tilde
 D,$ s.c.c. and with $\displaystyle {\mathcal{C}}^{3}$boundary
 such that we have $\displaystyle \ \tilde D\subset U$ and $\displaystyle
 \partial D\cap B(\zeta ,r)=\partial \tilde D\cap B(\zeta ,r).$
\end{Prps}
\quad Proof.\ \par 
Let $\rho $ be a defining function for $\displaystyle D.$ Let
 $\displaystyle \zeta \in \partial D$ and $U$ a neighbourhood
 of $\zeta $ in $\displaystyle {\mathbb{C}}^{n}.$ Consider a
 positive convex increasing function $\chi $ defined on $\displaystyle
 {\mathbb{R}}^{+},\ \ {\mathcal{C}}^{\infty }$ and such that
 $\displaystyle \chi =0$ in $\displaystyle (0,r)$\!\!\!\! . Set
 $\displaystyle \tilde \rho (z):=\rho (z)+a\chi (\left\vert{z-\zeta
 }\right\vert ^{2})$ ; we have $\displaystyle \partial \bar \partial
 \tilde \rho =\partial \bar \partial \rho +a\partial \bar \partial
 \chi .$ But, as is easily seen, $\displaystyle i\partial \bar
 \partial \chi $ is positive at each point $\displaystyle z,$
 hence, setting $\displaystyle A=i\partial \bar \partial \rho
 ,\ B=ai\partial \bar \partial \chi ,$ we can apply lemma~\ref{SE6}
 and we have that the domain $\displaystyle \tilde D:=\lbrace
 \tilde \rho <0\rbrace $ is also s.c.c. with smooth $\displaystyle
 {\mathcal{C}}^{3}$ boundary.\ \par 
Now we choose $\displaystyle r$ small enough to have $\displaystyle
 B(\zeta ,3r)\subset U.$ We have $\displaystyle \tilde \rho (z)<0\Rightarrow
 \rho (z)<-a\chi (\left\vert{z-\zeta }\right\vert ^{2})$ ; so we set :\ \par 
\quad \quad \quad $\displaystyle \alpha :=\sup \ _{z\in D}-\rho (z)<\infty ,$ by
 the compactness of $\bar D$ and $\displaystyle \beta :=\inf
 \ _{z\in U\backslash B(\zeta ,2r)}\chi (\left\vert{z-\zeta }\right\vert
 ^{2})=4r^{2}.$\ \par 
Then with $\displaystyle a:=\frac{\alpha +1}{\beta }$ we get
 that $\displaystyle \lbrace \tilde \rho (z)<0\rbrace \subset
 U$ because if not $\displaystyle \exists z\notin B(\zeta ,3r)::\rho
 (z)<-a\chi (\left\vert{z-\zeta }\right\vert ^{2})<-(\alpha +1)$
 which is not possible.\ \par 
\quad Of course in the ball $\displaystyle B(\zeta ,r)$ we have $\displaystyle
 \partial D\cap B(\zeta ,r)=\partial \tilde D\cap B(\zeta ,r).$
 $\displaystyle \blacksquare $\ \par 
\quad We shall need to extend this proposition to the case of $\displaystyle
 {\mathcal{C}}^{3}\ c$-convex intersection. \ \par 

\begin{Prps}
~\label{SI16}(Localizing s.c.c. intersection)Let $D$ be a $\displaystyle
 {\mathcal{C}}^{3}\ c$-convex intersection  in $\displaystyle
 {\mathbb{C}}^{n}.$ Let $\displaystyle \zeta _{0}\in \partial
 D,\ U$ a neighbourhood of $\zeta _{0}$ in $\displaystyle {\mathbb{C}}^{n}$
 and $\displaystyle B(\zeta _{0},r)$ a ball centered at $\zeta
 $ and of radius $r$ ; then there is a domain $\displaystyle
 \ \tilde D,\ {\mathcal{C}}^{3}\ c$-convex intersection, such
 that we have $\displaystyle \ \tilde D\subset U$ and $\displaystyle
 \partial D\cap B(\zeta _{0},r)=\partial \tilde D\cap B(\zeta _{0},r).$
\end{Prps}
\quad Proof.\ \par 
By assumption the $1$-forms $\displaystyle \lbrace d\rho _{j}(\zeta
 )\rbrace _{j\in I}$ are linearly independent in $\displaystyle
 \ \bigcap_{j\in I}{\lbrace \rho _{j}(z)\leq 0\rbrace }$ and
 $\displaystyle \ \left\vert{I}\right\vert \leq n-3.$ Take a
 point $\zeta _{0}\in \bigcap_{j\in I}{\lbrace \rho _{j}(z)=0\rbrace
 },$ by translation in $\displaystyle {\mathbb{C}}^{n},$ we may
 suppose that $\displaystyle \zeta _{0}=0,$ and we have to define
 the domain $\ \tilde D$ with the properties stated in the proposition.\ \par 
Take a vector $\displaystyle h\in {\mathbb{C}}^{n}$ of norm $1$
 and set, with $\displaystyle a\cdot b:=\sum_{k=1}^{n}{a_{k}b_{k}},$\ \par 
\quad \quad \quad $\displaystyle \rho (z):=(\left\vert{z}\right\vert ^{2}-r^{2})(1+h\cdot
 z+\bar h\cdot \bar z).$\ \par 
Because $\displaystyle 1+h\cdot z+\bar h\cdot \bar z>0$ near
 the origin, $\rho (z)$ is a defining function for the ball centered
 at $0$ and of radius $\displaystyle r.$\ \par 
\quad Now the claim is: we can choose the vector $\displaystyle h$
 in such a way that $\displaystyle d\rho _{j}(0)$ and $d\rho
 (0)$ are linearly independent for $\displaystyle j\in I.$\ \par 
\quad We have $\displaystyle d\rho (0)=(-r^{2})(h\cdot dz+\bar h\cdot
 d\bar z).$\ \par 
We already know that the $\displaystyle d\rho _{j}(0)$ are linearly
 independent and span a space of dimension less than $\displaystyle
 n-3,$ so we take $\displaystyle h$ in such a way that the form
 $\displaystyle (h\cdot dz+\bar h\cdot d\bar z)$ is not in the
 span of the $\displaystyle d\rho _{j}(0)$ for $\displaystyle
 j\in I.$ This is independent of the choice of $\displaystyle
 r>0.$ By continuity this is still true for $z$ in a neighbourhood
 $V$ of $0$ with $V$ independent of $\displaystyle r>0,$ so we
 now choose $\displaystyle r>0$ in order that $\displaystyle
 B(0,r)\subset V.$ We extend $\rho $ outside of the ball $\displaystyle
 B(0,r)$ to be a $\displaystyle {\mathcal{C}}^{\infty }$ function
 $\tilde \rho $ in $\displaystyle {\mathbb{C}}^{n},\ \tilde \rho
 =\rho $ in $\displaystyle B(0,r)$ and strictly positive in $\displaystyle
 \bar B(0,r)^{c}$ in order for this $\tilde \rho $ to be a genuine
 defining function for $\displaystyle B(0,r).$ So we get the
 condition (i) in the definition~\ref{SI14}.\ \par 
\quad The condition (ii) is easier because we have that\ \par 
\quad \quad \quad $\displaystyle \partial \bar \partial \tilde \rho =(1+h\cdot
 z+\bar h\cdot \bar z)\sum_{k=1}^{n}{dz_{k}\wedge d\bar z_{k}}+h\cdot
 dz\wedge z\cdot d\bar z+\bar z\cdot dz\wedge \bar h\cdot d\bar
 z=\sum_{k=1}^{n}{dz_{k}\wedge d\bar z_{k}}+{\mathcal{O}}(\left\vert{z}\right\vert
 ).$\ \par 
Hence there is a subspace $\displaystyle T_{z}^{I}$ such that
 for $\displaystyle i\in I$ the Levi forms $\displaystyle L\rho
 _{i}$ restricted on $\displaystyle T_{z}^{I}$ are positive definite
 by hypothesis and, because $\displaystyle L\tilde \rho $ is
 positive definite everywhere, we have that it is also positive
 definite on $\displaystyle T_{z}^{I}$ which has the right dimension
 $\displaystyle n-c+1.$\ \par 
It is at this point that we need $\displaystyle N\leq n-3,$ because
 we add the new domain $\displaystyle B(0,r).$ $\displaystyle
 \blacksquare $\ \par 
\ \par 
{\bf On a theorem of H. Rossi.}\ \par 
\ \par 
\quad We shall use the following lemma.\ \par 

\begin{Lmm}
~\label{SE5} Let $\displaystyle A,\ B$ two self adjoint $\displaystyle
 n{\times}n$ matrices such that $A$ has at least $\displaystyle
 d-c+1$ strictly positive eigenvalues and $\displaystyle \mathrm{k}\mathrm{e}\mathrm{r}A$
 is of dimension $\displaystyle n-d$ and $B$ is positive and
 has $\displaystyle n-d$ eigenvectors in $\displaystyle \mathrm{k}\mathrm{e}\mathrm{r}A$
 associated to strictly positive eigenvalues. Then $\displaystyle
 A+B$ has at least $\displaystyle n-c+1$ strictly positive eigenvalues.
\end{Lmm}
\quad Proof.\ \par 
Because $A$ is self adjoint, the spaces $\displaystyle \mathrm{k}\mathrm{e}\mathrm{r}A$
 and $\displaystyle H:=\mathrm{k}\mathrm{e}\mathrm{r}A^{\perp
 }$ are invariant for $\displaystyle A.$ Because $\displaystyle
 \mathrm{k}\mathrm{e}\mathrm{r}A$ has dimension $\displaystyle
 n-d$ and there is $\displaystyle n-d$ eigenvectors of $B$ in
 it, then $\displaystyle \mathrm{k}\mathrm{e}\mathrm{r}A$ is
 generated by these eigenvectors. Hence, because $B$ is self
 adjoint, this means that $\displaystyle \mathrm{k}\mathrm{e}\mathrm{r}A$
 and $H$ are also invariant for $B.$ Set $\displaystyle S:=A+B.$\ \par 
\quad Let $\displaystyle v\in \mathrm{k}\mathrm{e}\mathrm{r}A$ be such
 that $\displaystyle Bv=\lambda v,\ \lambda >0,$ then $\displaystyle
 Sv=Av+Bv=Bv=\lambda v$ ; hence on $\displaystyle \mathrm{k}\mathrm{e}\mathrm{r}A,\
 S$ has $\displaystyle n-d$ strictly positive eigenvalues.\ \par 
\quad On $H$ we have $\displaystyle B\geq 0$ and $A$ has at least $\displaystyle
 d-c+1$ strictly positive eigenvalues, hence on $H$ we can apply
 lemma~\ref{SE6} and we have that $S$ has at least $\displaystyle
 d-c+1$  strictly positive eigenvalues on $\displaystyle H.$
 Because $H$ and $\displaystyle \mathrm{k}\mathrm{e}\mathrm{r}A$
 have an intersection reduced to $\displaystyle \lbrace 0\rbrace
 ,\ S$ has $\displaystyle d-c+1+n-d=n-c+1$  strictly positive
 eigenvalues. $\displaystyle \blacksquare $\ \par 
\quad The aim is to extend a theorem by H. Rossi~\cite{Rossi76} where
 we replace strictly pseudo convex by strictly $c$-convex.\ \par 

\begin{Thrm}
~\label{SE9}Let $M$ be a closed submanifold of a Stein domain
 $\displaystyle U_{0}$ in $\displaystyle {\mathbb{C}}^{n}.$ Suppose
 there is a neighbourhood $U$ of $M$ and an holomorphic retraction
 $\displaystyle \pi :\ U\rightarrow M.$ Let $D$ be a strictly
 $c$-convex domain in $M,\ \bar D\subset M.$\par 
\quad Then there is a strictly $c$-convex domain $E$ in $\displaystyle
 {\mathbb{C}}^{n}$ such that :\par 
\quad (A)  $\displaystyle \bar E\subset U\cap U_{0}$\par 
\quad (B) $\displaystyle E\cap M=D$\par 
\quad (C) $\displaystyle \partial E$  cuts $M$ transversely along $\displaystyle
 \partial D$\par 
\quad (D) $\displaystyle \pi :\ \bar E\rightarrow \bar D.$
\end{Thrm}
\quad Proof.\ \par 
I shall copy the main points in the proof by H. Rossi making
 the necessary changes.\ \par 
\quad Docquier and Grauert (see ~\cite{Rossi76}) give us a neighbourhood
 $U$ of $\displaystyle \bar D$ in $\displaystyle {\mathbb{C}}^{n}$
 and a retraction $\displaystyle \pi \ :\ U\rightarrow M\cap
 U$ such that the fibers of $\pi $ cut transversely $\displaystyle
 M\cap U$ and are of dimension $\displaystyle n-d.$\ \par 
\quad We set for $\displaystyle z\in U$ and $\displaystyle j=1,\cdots
 ,n,\ f_{j}(z)=z_{j}-\pi _{j}(z).$ The equations $\displaystyle
 z-\pi (z)=0$ define the sub manifold $M$ :\ \par 
if $\displaystyle z\in M,\ \pi (z)=z$ because $\pi $ is a retraction
 {\sl on} $M$ ; if $\displaystyle z\notin M,$ because $\displaystyle
 \pi (z)\in M,z-\pi (z)\neq 0.$ Moreover, because the fibers
 of $\pi $ cut transversely $M$ at any point $\zeta $ of $\bar
 D,$ we have that the jacobian matrix contains a $\displaystyle
 (n-d){\times}(n-d)$ sub determinant which is not $0$ at $\zeta
 ,$ hence not $0$ in a neighbourhood of this point. This means
 that, by a change of variables, the set $\displaystyle (f_{j})_{j=1,\cdots
 ,n}$ contains a coordinates system for the fibers of $\pi $
 at any point of $\bar D,$ hence at all points of a neighbourhood
 $\displaystyle U_{1}$ of $\bar D$ in $\displaystyle {\mathbb{C}}^{n}.$
 These "explicit" functions replace the one generating the idealsheaf
 of $M$ used by H. Rossi.\ \par 
\quad Let $\rho $ be a defining function for $D$ in $M,$ we still follow
 H. Rossi and we set:\ \par 
\quad \quad \quad $\displaystyle \sigma (z):=\rho \circ \pi +A\sum_{j=1}^{n}{\left\vert{f_{j}}\right\vert
 ^{2}},$\ \par 
where the constant $A$ will be chosen later. Because $\displaystyle
 F(z):=\sum_{j=1}^{n}{\left\vert{f_{j}(z)}\right\vert ^{2}}=0$
 on $\displaystyle M\cap U,$ it exists a $\displaystyle \epsilon
 _{0}>0$ such that $\displaystyle \lbrace F(z)<\epsilon _{0}\rbrace
 \cap U\subset U_{1}.$\ \par 
\quad It remains to see that $\sigma $ is strictly $c$-convex, i.e.
 $\displaystyle i\partial \bar \partial \sigma $ has at least
 $\displaystyle n-c+1$  strictly positive eigenvalues.\ \par 
Fix $\displaystyle \zeta \in \bar D$ ; because $D$ is strictly
 $c$-convex, $\displaystyle i\partial \bar \partial \rho \circ
 \pi (\zeta )$ has at least $\displaystyle d-c+1$  strictly positive
 eigenvalues on the tangent space to $M$ at $\zeta .$ Because
 the set $\displaystyle (f_{j})_{j=1,\cdots ,n}$ contains a coordinates
 system for the fibers of $\pi $ we have $\displaystyle i\partial
 \bar \partial (\sum_{j=1}^{n}{\left\vert{f_{j}}\right\vert ^{2}})$
 has all, i.e. $\displaystyle n-d,$  strictly positive eigenvalues
 on the tangent space to the fiber $\displaystyle \pi ^{-1}\pi
 (\zeta )$ at $\displaystyle \zeta .$\ \par 
Because the kernel of $\displaystyle i\partial \bar \partial
 \rho \circ \pi $ is the tangent space to the fiber $\displaystyle
 \pi ^{-1}\pi (\zeta ),$ we get, by lemma~\ref{SE5}, that $\displaystyle
 i\partial \bar \partial \sigma =i\partial \bar \partial \rho
 \circ \pi +i\partial \bar \partial (\sum_{j=1}^{n}{\left\vert{f_{j}}\right\vert
 ^{2}})$ has at least $\displaystyle n-c+1$  strictly positive
 eigenvalues. So we have at least $\displaystyle n-c+1$  strictly
 positive eigenvalues at any point of $\bar D$ hence also in
 a neighbourhood $V$ of $\bar D$ in $\displaystyle {\mathbb{C}}^{n}.$
 Now we take $\displaystyle A\epsilon _{0}>\sup \ _{z\in D}\left\vert{\rho
 (z)}\right\vert $ and we set $\displaystyle E:=\lbrace z\in
 U\cap V::\sigma (z)<0\rbrace $ ; we get exactly as H. Rossi,
 that $E$ is strictly $c$-convex and we have all properties of
 the theorem. $\displaystyle \blacksquare $\ \par 
\quad We have to get an analogous result in the case where $D$ is a
 $\displaystyle {\mathcal{C}}^{3}\ c$-convex intersection  in
 $\displaystyle M.$ \ \par 

\begin{Thrm}
~\label{SI18}Let $M$ be a closed submanifold of a Stein domain
 $\displaystyle U_{0}$ in $\displaystyle {\mathbb{C}}^{n}.$ Suppose
 there is a neighbourhood $U$ of $M$ and an holomorphic retraction
 $\displaystyle \pi :\ U\rightarrow M.$ Let $D:=\bigcap_{k=1}^{N}{D_{k}}$
 be a $\displaystyle {\mathcal{C}}^{3}\ c$-convex intersection
 in $M,\ \bar D\subset M.$\par 
\quad Then there is a $\displaystyle {\mathcal{C}}^{3}\ c$-convex intersection
 $E:=\bigcap_{k=1}^{N}{\tilde D_{k}}$ in $\displaystyle {\mathbb{C}}^{n}$
 such that :\par 
\quad (A)  $\displaystyle \bar E\subset U\cap U_{0}$\par 
\quad (B) $\displaystyle E\cap M=D$\par 
\quad (C) $\displaystyle \partial E$  cuts $M$ transversely along $\displaystyle
 \partial D$\par 
\quad (D) $\displaystyle \pi :\ \bar E\rightarrow \bar D.$
\end{Thrm}
\quad Proof.\ \par 
By theorem~\ref{SE9}, and with the same notations, we can extend
 each $\displaystyle D_{k}$ by $\displaystyle \tilde D_{k}:=\lbrace
 \tilde \rho _{k}<0\rbrace $ in $\displaystyle {\mathbb{C}}^{n},$
 with $\displaystyle \tilde \rho _{k}:=\rho _{k}\circ \pi +AF,$
 where $\displaystyle \rho _{k}$ are the defining function for
 $\displaystyle D_{k}$ and $\displaystyle F(z):=\sum_{j=1}^{n}{\left\vert{f_{j}(z)}\right\vert
 ^{2}},$ such that they fulfil the $\displaystyle {\mathcal{C}}^{3}\
 c$-convex intersection requirements in $\displaystyle M.$ The
 point is to see that we can choose $A$ in such a way that $\tilde
 D:=\bigcap_{k=1}^{N}{\tilde D_{k}}$ fulfils the $\displaystyle
 {\mathcal{C}}^{3}\ c$-convex intersection requirements in $\displaystyle
 {\mathbb{C}}^{n}.$\ \par 
\quad First we choose $\displaystyle A\epsilon _{0}>\sup \ _{k=1,...,N,\
 z\in D_{k}}\left\vert{\rho _{k}(z)}\right\vert $ in order to
 have that all $\displaystyle \tilde D_{k}$ are in the domain
 of the retraction $\pi $ as for theorem~\ref{SE9}.\ \par 
\quad Fix a point $\displaystyle z_{0}\in \tilde D$ and take a vector
 $X$ in $\displaystyle {\mathbb{C}}^{n};$ then we can decompose
 it as $\displaystyle X=X_{M}\oplus X_{F}$ where $X_{M}$ is tangent
 at $z_{0}$ to the manifold $\displaystyle \lbrace z::z-\pi (z)=z_{0}-\pi
 (z_{0})\rbrace ,$ "parallel to $M$", and $\displaystyle X_{F}$
 is tangent to the fiber passing through $\displaystyle z_{0},\
 \lbrace z::\pi (z)=\pi (z_{0})\rbrace ,$ because we know that
 the fibers are transverse to $M,$ which is still true in a neighbourhood
 $V$ of $\bar D$ in $\displaystyle {\mathbb{C}}^{n}.$ Choose
 $A$ big enough to have all the $\tilde D_{k}$ in $\displaystyle
 V,$ the same way we did it in the proof of theorem~\ref{SE9}.\ \par 
\ \par 
Now if $\displaystyle z\in \bar D\subset M$ we already have that
 the $\displaystyle \lbrace d\tilde \rho _{j}(z)\rbrace _{j\in
 I}$ are linearly independent because there we have $\displaystyle
 F(z)=dF(z)=0$ hence $\displaystyle d\tilde \rho _{j}(z)=d\rho
 _{j}(z).$ So we make the assumption that $\displaystyle z\notin
 M.$ Let $\displaystyle I=(i_{1},...,i_{l})$ and suppose that
 the $\displaystyle \lbrace d\tilde \rho _{j}\rbrace _{j\in I}$
 are not linearly independent, then there is $\displaystyle \lambda
 \in {\mathbb{R}}^{\left\vert{I}\right\vert }$ such that\ \par 
\quad \quad \quad $\displaystyle \exists z\in \tilde D::0=\sum_{j\in I}{\lambda
 _{j}d\tilde \rho _{j}(z)}=\sum_{j\in I}{\lambda _{j}d(\rho _{j}\circ
 \pi (z)}+(\sum_{j\in I}{\lambda _{j}})AdF(z).$\ \par 
This means that\ \par 
\quad \quad \quad \begin{equation}  -(\sum_{j\in I}{\lambda _{j}})AdF(z)=\sum_{j\in
 I}{\lambda _{j}d(\rho _{j}\circ \pi (z)}.\label{SI22}\end{equation}\ \par 
Take any vector $\displaystyle X$ tangent to $\displaystyle {\mathbb{C}}^{n}$
 at $z$ ; then we have $\displaystyle X=X_{M}\oplus X_{F}$ and
 $\displaystyle \ {\left\langle{d(\rho _{j}\circ \pi )(z),\ X_{F}}\right\rangle}=0$
 because $\rho _{j}\circ \pi (\zeta )$ is constant along the
 fiber $\displaystyle \lbrace \zeta ::\pi (\zeta )=\pi (z)\rbrace
 .$ The same way $\displaystyle \ {\left\langle{dF(z),\ X_{M}}\right\rangle}=0$
 because $\displaystyle F(z)=\left\vert{z-\pi (z)}\right\vert
 ^{2}$ is constant along $\displaystyle \lbrace \zeta ::\zeta
 -\pi (\zeta )=z-\pi (z)\rbrace .$ So\ \par 
\quad \quad \quad $\displaystyle -(\sum_{j\in I}{\lambda _{j}})A{\left\langle{dF(z),\
 X}\right\rangle}=-(\sum_{j\in I}{\lambda _{j}})A{\left\langle{dF(z),\
 X_{F}}\right\rangle}=$\ \par 
\quad \quad \quad \quad \quad $\displaystyle =\sum_{j\in I}{\lambda _{j}{\left\langle{d(\rho
 _{j}\circ \pi )(z),\ X_{F}}\right\rangle}}=0.$\ \par 
Hence for any $\displaystyle X\in {\mathbb{C}}^{n},\ (\sum_{j\in
 I}{\lambda _{j}})A{\left\langle{dF(z),\ X}\right\rangle}=0$
 which means that the $1$-form $\displaystyle (\sum_{j\in I}{\lambda
 _{j}})AdF(z)=0$ hence, because $\displaystyle dF(z)\neq 0$ for
 $\displaystyle z\notin M,$ we have that $\displaystyle (\sum_{j\in
 I}{\lambda _{j}})=0.$\ \par 
\quad But by~(\ref{SI22}) this implies $\displaystyle \ \sum_{j\in
 I}{\lambda _{j}d(\rho _{j}\circ \pi (z)}=0$ which means that
 $\displaystyle \lambda _{j}=0$ because the $\displaystyle d\rho
 _{j}(\zeta )$ are independent at all points and in particular
 at the point $\displaystyle \zeta =\pi (z).$ So a contradiction
 which proves that the $\displaystyle d\tilde \rho _{k}$ are
 linearly independent.\ \par 
\quad To have the ii) fix $\displaystyle z\in \bigcap_{j\in I}{\lbrace
 \tilde \rho _{j}\leq 0\rbrace },$ and set $\displaystyle \zeta
 :=\pi (z)\in \bar D.$ The points $\displaystyle z,\zeta $ belongs
 to an open set $\displaystyle U:=U_{j}$ of the covering $\displaystyle
 (U_{j},\varphi _{j})$ done via lemma~\ref{SI20}, so reading
  by $\displaystyle \varphi :=\varphi _{j}$ we are in the following
 situation (I keep the same notations) : we have $\displaystyle
 z=(z',z"),\ \ \zeta =(z',0)$ and the retraction $\pi $ is the
 orthogonal projection $\displaystyle w\rightarrow (w',0)$ where
 $\displaystyle w':=(w_{1},...,w_{d})\ ;\ w":=(w_{d+1},...,\
 w_{n}).$ The tangent space $\displaystyle T_{\zeta }(M)$ is
 just $\displaystyle \lbrace w::w"=0\rbrace $ and by the hypotheses
 on the $\displaystyle \rho _{j}(w)=\rho _{j}(w')$ we know that
 there is a subspace $\displaystyle T_{\zeta }^{I},$ of dimension
 at least $\displaystyle d-c+1,$ of the tangent space $\displaystyle
 T_{\zeta }(M)$ on which the Levi forms $\displaystyle L\rho
 _{j}(\zeta )$ are positive definite. Lifting this space $\displaystyle
 T_{\zeta }(M)$ at the point $z$ keeping it parallel to itself,
 call it $\displaystyle T_{z}^{I},$  it still have dimension
 $\displaystyle d-c+1,$ and because the $\displaystyle \rho _{j}$
 do not depend on $\displaystyle w",$ we still have that the
 Levi form $\displaystyle L(\rho _{j}\circ \pi )(z)$ on $\displaystyle
 T_{z}^{I}$ is the same as the Levi form $\displaystyle L\rho
 _{j}(\zeta )$ on $\displaystyle T_{\zeta }^{I},$ so it is positive
 definite.\ \par 
\quad Now we have $\displaystyle \tilde \rho _{k}:=\rho _{k}\circ \pi
 +AF$ and $\displaystyle i\partial \bar \partial F$ has all its
 eigenvalues positive so on $\displaystyle T_{z}^{I}$ the Levi
 form $\displaystyle L\rho _{j}(z)$ is positive definite by the
 proof of lemma~\ref{SE6}. $\displaystyle \blacksquare $\ \par 
\ \par 

\bibliographystyle{/usr/local/texlive/2013/texmf-dist/bibtex/bst/base/plain}

\end{document}